\documentclass{article}

\usepackage{amssymb}

\usepackage[english,francais]{babel}

\newtheorem{theorem}{Theorem}[section]

\newtheorem{e-proposition}[theorem]{Proposition}

\newtheorem{e-definition}[theorem]{Definition\rm}

\newtheorem{theoreme}{Th\'eor\`eme}[section]

\newtheorem{proposition}[theoreme]{Proposition}

\newtheorem{exemple}{\it Exemple\/}

\def\og{\leavevmode\raise.3ex\hbox{$\scriptscriptstyle\langle\!\langle$~}}
\def\fg{\leavevmode\raise.3ex\hbox{~$\!\scriptscriptstyle\,\rangle\!\rangle$}}
\title{Godbillon-Vey classes for super-foliations }

\author{Camille Laurent-Gengoux \\ Department of Math,\\ Penn State University,\\ University Park,\\ PA 16802 USA \\ laurent@math.psu.edu }

\begin{document}



\selectlanguage{english}

\vspace{-2.6cm}

\selectlanguage{francais}


\selectlanguage{english}


\maketitle

\begin{abstract}

We construct the equivalent of the Godbillon-Vey class
and  its generalizations
for regular foliations on super-manifolds.
We interpret these classes as classes of foliated flat connections.


\vskip 0.5\baselineskip

\selectlanguage{francais}
\noindent{\bf R\'esum\'e}
\vskip 0.5\baselineskip
\noindent
Nous d\'eterminons  l'\'equivalent de
la classe de Godbillon-Vey et de ses g\'en\'eralisations
pour des feuilletages r\'eguliers 
sur des super-vari\'et\'es.
Nous int\'erpr\`etons ensuite ces classes
comme des classes caract\'eristiques associ\'ees
\`a des connections feuillet\'ees plates.


\end{abstract}



La classe de Godbillon-Vey est un el\'ement de $H^3(M,{\mathbb R})$
associ\'ee \`a un feuilletage r\'egulier de codimension~$1$
sur une vari\'et\'e $M$. Depuis les travaux de Fuks, Bershtein ou
Rozenfel'd \cite{Bernshtein,Fuks2}, on sait plus g\'en\'eralement,
  associer \`a tout feuilletage de codimension $n$
sur une vari\'et\'e $M$ un homomorphisme
de $H^*\big(Vect(n)\big)$ \`a valeurs dans $H^*(M,{\mathbb R}) $,
ou, en d'autres mots, on sait construire autant de classes caract\'eristiques
qu'il y a de classes de cohomologie dans la cohomologie
de Gelfand-Fuks de l'alg\`ebre de Lie $Vect(n)$
des champs de vecteurs formels.

Consid\'erons ${\mathcal M}$ une super-vari\'et\'e
et $M$ sa vari\'et\'e de base.
Nous associons dans la pr\'esente note un \'el\'ement
de $H^1(M,{\mathbb R})$ \`a des super-feuilletages ${\mathcal F}$ de codimension
$0 + \epsilon 1$ par le th\'eor\`eme suivant:
\vspace{0.2cm}

{\bf Th\'eor\`eme A} {\em  Soit un super-feuillatege
${\mathcal F}$ de codimension $0+ \epsilon 1 $ d\'etermin\'e  
par une $1$-forme paire $a$ et soit ${\mathcal X}$ un super-champ
de vecteur tel que $a({\mathcal X})=1$. La classe de cohomologie
de $p(L_{\mathcal X}a) \in \Omega^1(M)$, o\`u $p:\Omega({\mathcal M})\to \Omega(M)$ est la projection canonique, ne d\'epend pas des choix de $a$ et ${\mathcal X}$.   }
\vspace{0.2cm}

 C'est cette classe que disons \^etre un \'equivalent
de la classe de Godbillon-Vey.
Plus g\'en\'eralement, nous affirmons que l'on peut
associer \`a un super-feuilletage au compl\'ementaire trivial
une application de $H^*\big( Vect(n,m)_0\big)$
dans $H^*(M,{\mathbb R}$, o\`u $H^*\big( Vect(n,m)_0\big)$
est la cohomologie de Gelfand-Fuks de la partie paire
$Vect(n,m)_0$ de la super-alg\`ebre de Lie $Vect(n,m)$
des super-champs de vecteurs formels \`a $n$ param\`etres pairs
et $m$ param\`etres impairs.
En d'autre mot, \'etant donn\'e un  $C \in H^*\big( Vect(n,m)_0\big)$,
nous construisons une classe caract\'eristique $\phi_{\mathcal F}(C) \in H^*(M,{\mathbb R})$.
\vspace{0.2cm}

{\bf Th\'eor\`eme B} {\em  Soit un super-feuillatege
${\mathcal F}$ de codimension $n+ \epsilon m $ admettant
un suppl\'ementaire trivialisable. 
Il existe un homomorphisme $\phi_{\mathcal F} :H^*\big( Vect(n,m)_0
\big) \to H^*(M)$ tel que
\begin{enumerate}
\item si $f^*:{\mathcal N}\to {\mathcal M} $ est une submersion,
alors $ \phi_{f^*({\mathcal F})}=\tilde{f} \circ  \phi_{\mathcal F}$,
o\`u $ f^*({\mathcal F})$ est le tir\'e en arri\`ere du super-feuilletage ${\mathcal F}$ et o\`u $\tilde{f}:N \to M $ est l'application lisse induite
par la restriction de $f$ aux vari\'et\'es de bases,
\item si ${\mathcal M}$ est une vari\'et\'e diff\'erentielle ordinaire,
 $\phi_{\mathcal F}$ est l'application de Bershtein, Fuks et Rozenfel'd,
 \cite{Bernshtein,Fuks2}.
\end{enumerate} }
\vspace{0.2cm}

Il peut sembler surprenant que ce soit la cohomologie $H^*(M)$ 
de la vari\'et\'e de base qui joue un r\^ole ici: cela est d\^u \`a un r\'esultat de Batchelor
\cite{Leites}
selon lequel cette cohomologie est isomorphe \`a celle de la super-vari\'et\'e ${\mathcal M}$.
Enfin, remarquons que c'est la cohomologie de la partie
paire $Vect(n,m)_0$ et non celle de $Vect(n)$
qui g\'en\'eralise le r\'esultat classique. Il est en fait heureux qu'il
en soit ainsi puisque l'on peut montrer \cite{Astiachkevitch,Fuks} que lorsque $n<m$
l'espace vectoriel $H^*\big(Vect(n,m)_0\big)$ n'est que de dimension $2$.

Nous expliquerons enfin pourquoi \`a partir de
 tout super-feuilletage sur une super-vari\'et\'e ${\mathcal M}$,
on peut construire
 ce que l'on appelle des connections feuillet\'ees
plates sur la vari\'et\'e de base $M$;
 c'est \`a dire un feuilletage sur la vari\'et\'e de base,
un fibr\'e vectoriel et une connection plate de 
ce fibr\'e vectoriel d\'efinie uniquement
le long des feuilles du feuilletage.
On obtient alors
\vspace{0.2cm}

{\bf Th\'eor\`eme C} {\em Si \`a deux super-feuilletages ${\mathcal F}$
et ${\mathcal F}'$ sont associ\'ees des connnections feuillet\'ees plates
isomorphes, alors $\phi_{\mathcal F}=\phi_{{\mathcal F}'} $.} 

\vspace{0.2cm}
 
En d'autres termes, nos classes caract\'eristiques peuvent
se comprendre comme des classes caract\'eristiques de connections feuillet\'ees plates.

\vspace{0.2cm}

\selectlanguage{english}
\section{Introduction}

 The Godbillon-Vey class is an element of  $H^{3}(M,{\mathbb R})$ associated to a regular foliation of
codimension $1$ on a manifold $M$. More generally, to any regular foliation $F$
with a trivial supplement 
on a manifold $M$ is uniquely associated a (non-trivial in general) homomorphism from $H^*\big( Vect(n) \big)$ to $H^*(M,{\mathbb R}) $,
 see \cite{Bernshtein} or  \cite{Fuks2}
(where $H^*\big( Vect(n)\big) $ is the 
Gelfand-Fuks
cohomology of the Lie algebra $Vect(n)$ of formal vector fields with $n$ parameters).

Let $ {\mathcal M}$ be a super-manifold and $M$ its basis manifold.
We construct in this article 
an equivalent of the Godbillon-Vey class on a super-manifold
${\mathcal M}$ endowed with a super-foliation of codimension $0+ \epsilon 1$.
This class is an element of $H^1(M,{\mathbb R})$. 
More generally,  we will show  how to associate  
to any super-foliation ${\mathcal F}$  of codimension $n+\epsilon m$ with  a trivial supplement,
 a homomorphism $\phi_{\mathcal F}: H^*\big( Vect(n,m)_0 \big) \to H^*(M,{\mathbb R})$ where
$H^*\big( Vect(n,m)_0 \big)$ is
 the Gelfand-Fuks cohomology 
of the even part $Vect(n,m)_0 $ of the super-Lie algebra 
$Vect(n,m)$ of formal vector fields with $n$ even parameters and $m$ odd parameters.
In other words, we define for any $ C \in H^*\big( Vect(n,m)_0 \big) $
a characteristic class $\phi_{\mathcal F}(C) \in H^*(M,{\mathbb R})$ of super-foliation with a trivial supplement.
 The reader must not be surprised by the fact that the de Rham
cohomology of $ {\mathcal M}$ is not used:
 by a result of Batchelor \cite{Leites}, this
cohomology is indeed isomorphic to the cohomology of the basis manifold. Note that only
the even part  $Vect(n,m)_0 $ appears: it is indeed better like this, since
 $H^*\big(Vect(n,m)\big)$  has only two non-trivial
generators for $n < m$, see \cite{Fuks2}.

We will explain then why these classes are in fact completely determined
by a foliation with a trivial supplement
on the basis manifold $M$, a trivial vector bundle and a flat foliated  connection
canonically associated to the initial super-foliation on ${\mathcal M}$
and can be therefore understood as characteristic classes of such objects.

\section{Foliations on super-manifolds}


 We denote by ${\mathcal O}({\mathcal M})$  the super-functions on a super-manifold  ${\mathcal M}$, by
$Vect({\mathcal M})$ the super-Lie algebra of super-vector fields and by
$ \big( \Omega({\mathcal M}),d \big)$ the super-algebra of differential forms
on $ {\mathcal M}$ endowed with its de Rham differential.


On the super-manifold ${\mathbb R}^{p,q}$, we denote by ${\mathcal B}$
the ${\mathcal O}({\mathbb R}^{p,q}) $-sub-module of $Vect({\mathbb R}^{p,q}) $ generated by the 
super-vector fields $(\frac{\partial}{\partial x_1},\dots,
 \frac{\partial}{\partial x_{p-n}}) $
and $(\frac{\partial}{\partial\theta_1},\dots,\frac{\partial}{\partial \theta_{q-m}}) $ where $x_1,\dots,x_p,\theta_1,\dots,\theta_q$
are the natural parametrisation of ${\mathbb R}^{p,q}$.
A {\it regular super-foliation ${\mathcal F}$  of codimension $n+\epsilon m$} is a sub-super-Lie algebra
        $Vect({\mathcal F})$ of $Vect({\mathcal M})$
such that for any $x$ in the basis manifold $ M$ there exists a neighborhood ${\mathcal U}$ of $x$ and a diffeomorphism
from  ${\mathcal U}$ and ${\mathbb R}^{p,q}$ that provides an isomorphism
between ${\mathcal B}$ and the restriction of  $Vect({\mathcal F})$  to ${\mathcal U}$.

We say that a super-foliation ${\mathcal F}$ admits a trivial supplement
when there exists  odd 
$1$-forms $\omega_1,\dots,\omega_n$ and even $1$-forms
$a_1,\dots,a_m $ generating a free ${\mathcal O}({\mathcal M})$-module
of super-dimension $ n + \epsilon m$
such that  
$Vect({\mathcal F})=\{{\mathcal X} \in Vect({\mathcal M})| \omega_1({\mathcal X})=\dots=
\omega_n({\mathcal X})=a_1({\mathcal X})=\dots=a_m{\mathcal X} =0 \}$.



An important point is that to any super-foliation ${\mathcal F}$
of codimension $n +\epsilon m$ on a super-manifold ${\mathcal M}$ is associated
by a canonical way a super-foliation with  a trivial supplement on a
$Gl(n,m)$-principal bundle over ${\mathcal M}$, see \cite{laurent}.
This $Gl(n,m)$-principal bundle is the super-manifold
of free $ n + \epsilon m$ families  $T^*(M) $ that vanishes on $ Vect({\mathcal F})$.
Therefore, we can always assume that the super-foliation admits
a trivial supplement, by replacing the initial
 super-foliation by this new one if necessary.

\section{An equivalent of Godbillon-Vey classes for super-foliations of codimension $0 + \epsilon 1$}

Let ${\mathcal F}$ be a regular super-foliation of codimension $ 0 + \epsilon 1$ with a trivial
supplement  given
by an even $1$-form $a$ and $\Theta$  an unique odd super-vector field satisfying $a(\Theta)=1$.
Let us consider the $1$-form on $M$ given by $ A = p(L_{\Theta}a)$
where $p: \Omega({\mathcal M}) \to  \Omega(M)$
is the natural projection.

\begin{theorem}\label{th1} The $1$-form $A$ is closed.
 Its cohomology class in $H^1(M,{\mathbb R})$ does not depend upon the choice of $a$ and $\Theta$.
\end{theorem}

This Theorem allows us to consider $[A] \in H^1(M,{\mathbb R})$ as a characteristic class
of a super-foliation of codimension $ 0 + \epsilon 1$.

{\bf Proof:}
Since the $1$-form $a$ defines a super-foliation of codimension $ 0+\epsilon 1$,
there exists an odd $1$-form  $b$ such that $da=b \land a$, see \cite{Monterde}.
Since $b$ is odd, we have $b \wedge b=0 $ and the identity
   $0=d^2a= a\land b \land b + a \land db =  a \land db$ holds.
This implies that $db=0$ because,  $a$ being a regular even $1$-form, 
we have the surprising property that  $ a \wedge c  =0$ implies $ c=0$
for any $c \in \Omega({\mathcal M})$.

Now from $L_{\Theta} a =\mbox{\i}_{\Theta} da =  \mbox{\i}_{\Theta} a \wedge b = b + a \wedge \mbox{\i}_{\Theta} b  $,
we deduce immediately that $A=p(L_{\Theta} a )=p(b) $.
The $1$-form $A$ is thus a closed $1$-form on $M$ that does not depend on the choice of $\Theta$.

To show that $[A]$ does not depend on the choice of the $1$-form
$a$ chosen, we just have to see that if we replace $a$ by $f a$ for $f$ an even invertible 
 super-function on ${\mathcal M}$, we have $d(fa)=c \land f a $ with $c= b+\frac{df}{f}  $ and, as a consequence,

  $$p(c) = p(b) + p (\frac{df}{f})=p(b)+\frac{dp( f)}{p(f)} .$$

The super-function $f$ being invertible, $ p(f)$ never vanishes, so
  the $1$-form $p (\frac{df}{f})=d\big( ln(|p(f)|) \big)$ is  exact on $M$. $\triangleright$

\begin{exemple}\label{ex:super-cercles}.
On the super-circle $S^{1,1}$, let us denote $x\in S^1$ the even parameter and $\theta$ 
the odd parameter. For any $t \in {\mathbb R}$, let us consider the super-foliation of codimension $ 0 + \epsilon 1$ 
defined by the $1$-form $d\theta +t \theta dx$. Let us choose $\Theta=\frac{d}{d\theta}$.
We have $A=p(L_{\Theta}a)=t dx$ and $[A]=t \in H^1(S^1,{\mathbb R}) \simeq {\mathbb R}$. 
The characteristic class that we obtain 
in this case is thus non-zero if $t \neq 0$.
\end{exemple}

\section{Generalization}

The Lie super-algebra $Vect(n,m)$ of super-vector fields on ${\mathbb R}^{n,m}$
with polynomial coefficients is defined  for example in \cite{Astiachkevitch}: it is the Lie super-algebra
 of super-derivations of
${\mathbb R}[x_1,\dots,x_n] \otimes \bigwedge {\mathbb R}^m$
 where ${\mathbb R}[x_1,\dots,x_n]$ is the 
algebra polynomial functions on ${\mathbb R}^n$.
We denote by $H^*\big( Vect(n,m)_0 \big)$ the cohomology
of the even part of $Vect(n,m)$.
The following Theorem defines
{\em characteristic classes of a super-foliation with a trivial supplement}.

\begin{theorem}\label{th2} To any super-foliation
with a trivial supplement of codimension $n + \epsilon m $,
we can associate a natural and non-trivial in general 
 homomorphism $\phi_{\mathcal F}$  from  $H^*\big( Vect(n,m)_0 \big)$
to $H^*(M,{\mathbb R})$ such that
\begin{enumerate}
\item if  $f: {\mathcal N} \to {\mathcal M}$ is a submersion of super-manifolds, then
  $ \phi_{f^* {\mathcal F}} = \tilde{f}^* \phi_{ \mathcal F}$ where the super-foliation
$f^* {\mathcal F}$ is the pull-back of ${\mathcal F}$ by $f$
and $\tilde{f}:N\to M $ is the smooth map induced by $f$
on basis manifolds.
\item if ${\mathcal M}$ is an (ordinary) smooth manifold, then $m=0$ and $\phi_{\mathcal F}$
reduces to the well-known map of Fuks, Bershtein and Rozenfel'd \cite{Bernshtein}, \cite{Fuks}.
\end{enumerate}
\end{theorem}

The proof of this statement turns to be far away from being easy and will be published in \cite{laurent}.

\begin{exemple}. It is not difficult to see (see \cite{laurent}) that $ H^1\big( Vect(0,1)_0 \big)  \simeq {\mathbb R}$ and $ H^i\big( Vect(0,1) \big)=0 $
for $i \geq 2$. The characteristic class associated to $H^1\big( Vect(0,1) \big)$ is the class described in
Theorem \ref{th1}.
\end{exemple}

\section{Super-foliation on a super-manifold and flat foliated connection
on its basis manifold}

\subsection{Characteristic classes of trivial flat foliated connections}\label{se:ffc}

 Let $F$ be a foliation of codimension $n$ on
a smooth manifold $M$ and $E$ a vector bundle over $M$.
We recall that a {\it foliated connection}
is a bilinear map   $\nabla: Vect(F) \otimes \Gamma(E) \to \Gamma(E)$ 
that satisfies the usual properties of a connection. We say that 
this foliated connection is {\it flat} if $\forall X,Y \in Vect(F),
 \forall s \in \Gamma (E)$
we have  $ \nabla_X \nabla_Y s -\nabla_Y \nabla_X s-
\nabla_{[X,Y]} s=0 $.
By a {\em trivial} foliated connection $(F,E,\nabla)$, we mean a foliation $F$  with
a trivial supplement on $M$, a trivial vector bundle
$E$ over $M$ and  a foliated connection $\nabla$.

 There is a natural application from
the set of flat trivial foliated connections $(F,E,\nabla)$ 
to the set of  super-foliations ${\mathcal F}$  with a trivial supplement:
given a trivial foliated connection $(F,E,\nabla)$$ $ we define first 
the super-manifold ${\mathcal M}$ whose super-functions
 are the sections of $ \wedge^* E$
and the super-foliation generated by 
the  even derivations  $ \nabla_X$  of $ {\mathcal O}({\mathcal M})$ for all $X \in Vect(F)$.
By Theorem \ref{th2}, we can therefore
associate to any flat trivial foliated connection
a homomorphism from  $H^*\big( Vect(n,m)_0 \big)$ to $H^*(M,{\mathbb R})$.
This  defines 
{\em characteristic classes of a trivial flat foliated connection},
which can be shown to be non-trivial in general.

\subsection{A functor from super-foliations to flat foliated connections}

Given a  super-foliations with  a trivial supplement on a  super-manifold ${\mathcal M}$
one can construct a 
  flat trivial foliated connection as follows.
 Let  ${\mathcal F}$  be a super-foliation ${\mathcal F}$ 
 with  a trivial supplement.

\begin{itemize} 
\item  First  a regular foliation $F$ on $M$ is defined by the Lie algebra $p\big( Vect({\mathcal F})_0 \big) $, $p$ being here the natural projection
from the even part of $Vect({\mathcal M})$ onto $Vect(M)$ and $Vect({\mathcal F})_0$
being the even part of $Vect({\mathcal F})$. 
 
\item Second a vector-bundle $E_F$ is defined to be the vector bundle whose sections
are of the form $ \frac{{\mathcal J}_F  {\mathcal O}({\mathcal M})}{{\mathcal I}^2} $, where $ {\mathcal I}$
is the ideal of nilpotent elements of $ {\mathcal O}({\mathcal M})$
and $ {\mathcal J}_F$ the algebra of odd super-functions constant on the leaves 
of $ {\mathcal F}$ ({\it i.e. } the algebra of super-functions $f$ such that $ df (\tilde{X} ) =0$
for any $ \tilde{X}  \in Vect({\mathcal F})_0$).

\item The foliated connection is then defined $\forall X \in Vect(F) , \forall s \in \Gamma (E)$ by
$$  \nabla_X s = df (\tilde{X})  \hspace{1cm}  \mbox{mod}  \,{\mathcal I}^2  $$
where $f \in {\mathcal J}_F  {\mathcal O}({\mathcal M})$ is a super-function 
with $ s = f \, \, \, \mbox{mod}  \,{\mathcal I}^2$ and $ \tilde{X} \in Vect({\mathcal M})$
is an even super-vector field that projects naturally on $X$. We let the reader check that 
this  defines a foliated connection. 

\item If the super-foliation ${\mathcal F}$ has a trivial supplement then 
  the foliation $F$ has a trivial supplement also and  the vector bundle $E_F$ is trivial.

\end{itemize}

We have therefore constructed a trivial foliated connection.
 It is easy now to check that:

\begin{proposition}\cite{laurent}
The foliated connection $\nabla$ on $E_F$ is flat.
\end{proposition}

This last Proposition achieves our construction.
Note that this construction is not the inverse of the construction of 
section \ref{se:ffc}, see \cite{laurent}.

\begin{theorem} \label{th:second} \cite{laurent}
 Two super-foliations ${\mathcal F}$ and ${\mathcal F}'$ 
with a trivial supplement on a super-manifold ${\mathcal M}$ defining isomorphic
 flat trivial foliated connections have the same characteristic classes, {\em i.e.},
the homomorphisms $\phi_{\mathcal F}$ and  $\phi_{{\mathcal F}'}$ from $H^*\big( Vect(n,m)_0 \big)$
to $H^*(M,{\mathbb R})$ associated to these super-foliations are equal.
\end{theorem}


 The foliation $F$ previously defined has a trivial supplement:
therefore the theory of characteristic classes of \cite{Bernshtein,Fuks}
applies and provides a homomorphism $\phi_{F}:H^*\big(Vect(n)^*) \to H^*(M,{\mathbb R})$.
These characteristic classes are among those constructed in Theorem \ref{th2}.
More precisely, there is a natural inclusion $i:H^*\big(Vect(n)^*) \to H^*\big(Vect(n,m)_0^*)  $,
see \cite{laurent}, and  it can be shown that

\begin{proposition} \cite{laurent}
The following diagram is commutative
$$\begin{array}{ccc}
  H^*\big(Vect(n)^*)  &\stackrel{i}{\to}   &  H^*\big(Vect(n,m)_0^*) \\
    &\searrow   \phi_{F}& \downarrow \phi_{\mathcal F}\\
   &   & H^*(M,{\mathbb R}) \\
\end{array} . $$
\end{proposition}





\end{document}